%% file: template.tex
\documentclass{article}

\usepackage{arxiv}
\input{parameters}

\usepackage[utf8]{inputenc} 
\usepackage[T1]{fontenc}    
\usepackage{hyperref}       
\usepackage{url}            
\usepackage{booktabs}       
\usepackage{amsfonts}       
\usepackage{nicefrac}       
\usepackage{microtype}      
\usepackage{lipsum}
\usepackage{graphicx}
\usepackage{xcolor}
\graphicspath{ {./images/} }
\usepackage{amsmath}
\theoremstyle{definition}
\newtheorem{theorem}{Theorem}[section]

\newtheorem{definition}[theorem]{Definition}

\newtheorem{propriete}[theorem]{Propriety}

\newtheorem{lemma}[theorem]{Lemma}

\newtheorem{problem}[theorem]{Problem}

\newtheorem{proposition}[theorem]{Proposition}
\newtheorem{remark}[theorem]{Remark}

\title{A new geometric approach to multiobjective linear programming problems}
\author{
 Mustapha Kaci$^{{\color{white}.}\bf 1}$ and Sonia Radjef$^{{\color{white}.}\bf 2}$\\
 $^{\bf 1,2{\color{white}.}}$Department of Mathematics, Signal Image Parole (SIMPA) Laboratory\\ University of Oran Mohamed Boudiaf USTO-MB, Oran, Algeria.\\
  {\color{blue}\texttt{kaci.mustapha.95@gmail.com} }$^{{\color{white}.}\bf 1}$\\
	{\color{blue}\texttt{soniaradjef@yahoo.fr} }$^{{\color{white}.}\bf 2}$
   }

\begin{document}
\fontfamily{ptm}\selectfont
\maketitle
\begin{abstract}
 In this paper, we present a novel method for solving multiobjective linear programming problems (MOLPP) that overcomes the need to calculate the optimal value of each objective function. This method is a follow-up to our previous work on sensitivity analysis, where we developed a new geometric approach. The first step of our approach is to divide the space of linear forms into a finite number of sets based on a fixed convex polygonal subset of $\mathbb{R}^{2}$. This is done using an equivalence relationship, which ensures that all the elements from a given equivalence class have the same optimal solution. 
 We then characterize the equivalence classes of the quotient set using a geometric approach to sensitivity analysis. This step is crucial in identifying the ideal solution to the MOLPP. By using this approach, we can determine whether a given MOLPP has an ideal solution without the need to calculate the optimal value of each objective function. This is a significant improvement over existing methods, as it significantly reduces the computational complexity and time required to solve MOLPP.\\ 
 To illustrate our method, we provide a numerical example that demonstrates its effectiveness. Our method is simple, yet powerful, and can be easily applied to a wide range of MOLPP. This paper contributes to the field of optimization by presenting a new approach to solving MOLPP that is efficient, effective, and easy to implement.
\end{abstract}

\section{Introduction.}Linear programming is a widely used mathematical optimization technique in operational research and mathematical programming. It involves optimizing a mathematical program where the objective function and the functions defining the constraints are linear \cite{ref26,kaci2,ref2}. Linear constraints form a convex polyhedron, and the corners of the polyhedron are the basic feasible solutions, one of which can be the optimal solution. The results of convexity have been utilized to develop new numerical methods for resource allocation \cite{ref5}. In practical applications, linear programming is often used to model problems such as maximizing a company's profits, subject to various constraints such as resource limitations and market conditions. However, changes in market data can require updates to the initial model coefficients, making sensitivity analysis an essential part of linear programming, also known as post-optimal analysis \cite{ref1,ref29,ref3}. Sensitivity analysis is used to illustrate the range of linear program parameters for which the solution of the initial problem remains stable \cite{ref21,ref3}.
Multiobjective programming problems are mathematical problems in which one decision-maker seeks to optimize several generally conflicting objectives. Such problems arise in a variety of fields, including engineering, economics, and finance, where multiple criteria need to be considered simultaneously. Since the criteria space is typically supplied with a partial order, it is necessary to specify a sense of resolution (non-dominated solution) to develop solution methods \cite{ref28}. Multiobjective optimization techniques have been developed to enable an accurate analysis of trade-offs between competing objectives and to assist the decision-maker in reaching an acceptable trade-off \cite{ref2022}.\\
But nothing prevents all the objectives of a multiobjective problem from agreeing (achieving optimal value in an attainable common solution). To do so, we will develop an approach to recognize multiobjective issues that have no conflict (those that have an ideal solution). Then classify them as their solution, the classification will depend on the feasible region.
In this study, our goal is to address the problem of identifying whether a given MOLPP has an ideal solution where all of the objective functions are optimized simultaneously. To achieve this, we propose a new approach based on sensitivity analysis, building on the method presented in \cite{kaci}. Our method involves defining an equivalence relationship over the space of linear forms, which partitions the space into a finite number of equivalence classes. Each class contains linear forms that achieve their maximum value at a common point. This classification enables us to determine whether the MOLPP under consideration admits an ideal solution. Moreover, we introduce a new result that allows us to classify the MOLPP based on a fixed feasible region. 
In order to provide a comprehensive understanding of our approach, this document is structured into several sections.
Section \ref{Section_2} provides introductory remarks and presents the motivation for our work. In this section, we also provide an overview of the relevant literature and discuss related works that have investigated sensitivity analysis and multiobjective linear programming problems.
In Section \ref{Section_3}, we present the mathematical formulation of the MOLPP and discuss its properties, including the optimality of feasible solutions.
Section \ref{Section_4} outlines our proposed approach and presents a new result that allows us to classify the MOLPP on a fixed feasible region.
In Section \ref{Section_5}, we provide a numerical example to illustrate our method. Finally, in Section \ref{Section_6}, we summarize our findings and discuss potential directions for future research.

\section{\label{Section_2}Preliminaries.} This section gives a brief overview of our 2022 published work on the sensitivity analysis approach. For further details, please refer to \cite{kaci}. Our approach involves examining the impact of parameter or variable changes within a specified range on the optimal solution of a mathematical program. Our work builds upon existing research and introduces new techniques and results that we believe will contribute to the advancement of sensitivity analysis in multiobjective linear programming problems.
\par
Consider the linear programming problem in the standard form given below as our initial problem:
\begin{equation}\label{Pbm0}
\left\{
\begin{array}{l}
      \underset{x}\max\hspace{0.15cm} f(x)=c_{0}^{T}x           \\
      Ax\leq b\\
x\geq0
\end{array}.
\right.
\end{equation}
Where
$$
\begin{array}{cccc}
A=\left(\begin{array}{ccccccccc}
a_{11}&a_{12}\\
a_{21}&a_{22}\\
 \vdots    & \vdots \\
a_{m1}&a_{m2}
\end{array}\right),&b=\left(\begin{array}{c}b_{1}\\b_{2}\\\vdots\\b_{m}\end{array}\right),&c_{0}=\left(\begin{array}{c}c^{0}_{1}\vspace{0.15cm}\\c^{0}_{2}\end{array}\right),&x=\left(\begin{array}{c}x_{1}\\x_{2}\end{array}\right).
\end{array}
$$
and $f : \mathbb{R}^{2}\longmapsto\mathbb{R}$ is a linear form, $c_{0}$ is the vector of the objective function $f$.  $c^{0}_{1},c^{0}_{2}$ are constants, $x$ is a $2\times1$ vector of decision variables, $A$ is a $m\times 2$ matrix of constants,  $b$ is a $m\times 1$ vector of constants  and $m$ is the number of linear constraints.

\begin{problem}\label{Pbm1}

Let $x^{0}=(x^{0}_{1},x^{0}_{2})^{T}$ be an optimal solution of the problem (\ref{Pbm0}). A sensitivity analysis problem is to find all linear forms $g$ different from $f$  verifying:
\begin{equation}
\arg  \underset{x\in S} \max \, f(x)= \arg \underset{x\in S} \max \, g(x),\hspace{0.2cm} g\in\mathcal{L}\left(\mathbb{R}^{2}\right).
\end{equation}
Such that
$$
\arg  \underset{x\in S} \max \, f(x)=\left\{y\in S\hspace{0.2cm}|\hspace{0.2cm}f(y)\geq f(x),\hspace{0.2cm} x\in S\right\}.
$$
\end{problem}
Let us consider the two-dimensional Euclidean space $\mathbb{R}^{2}$ as the plane $z=0$. On this plane, we focus on the intersection of the graph of $f$ with the plane $z=0$, which is a two-dimensional subset of $\mathbb{R}^{2}$. We denote this subset as the line $d_{0}$, which can be represented as the set of all points $(x_{1},x_{2})$ in $\mathbb{R}^{2}$ for which $f(x_{1},x_{2}) = 0$.
\par
Geometrically, $d_{0}$ represents a linear subspace of $\mathbb{R}^{2}$ of dimension 1, since it is a line. This line is directed by the vector $v_{0}=(-c^{0}_{2},c^{0}_{1})$ that represents the slope of $d_{0}$ in the plane $z=0$. 

\begin{proposition}\label{proposition1}
Let $c_{1}^{0}, c_{2}^{0} \geq0$, then the problem (\ref{Pbm0}) is equivalent to:
\begin{equation}\label{Pbm2}
\underset{x\in S}\max \, \,\left\|x-P_{(d_{0})}(x)\right\|.
\end{equation}

\end{proposition}

\begin{problem}\label{Pbm3}
Consider the set $S$ defined as the convex hull of the extreme points $x^1$, $x^0$, and $x^2$, where $x^0$ is the optimal solution of problem (\ref{Pbm2}). We aim to find all vector lines $(d)$ that satisfy the following two inequalities simultaneously:
\begin{equation}\label{Pbm31}
                \left\|x^{0}-P_{(d)}\left(x^{0}\right)\right\|\geq\left\|x^{1}-P_{(d)}\left(x^{1}\right)\right\|,
\end{equation}
and
\begin{equation}\label{Pbm32}													
                \left\|x^{0}-P_{(d)}\left(x^{0}\right)\right\|\geq\left\|x^{2}-P_{(d)}\left(x^{2}\right)\right\| .	
   \end{equation}											
\end{problem}
\begin{proposition}
The problems (\ref{Pbm1}) and (\ref{Pbm3}) are equivalents.
\end{proposition}
\begin{remark}\label{remarktra} It is worth noting that the problem (\ref{Pbm3}) is constructed from another problem that has a solution. Therefore, it follows that the vector line $(d_{0})$ serves as a solution for the problem (\ref{Pbm3}), which can be considered a trivial solution. This is a significant observation, as it provides a starting point for exploring more complex solutions and identifying additional vector lines that satisfy the problem constraints.
\end{remark}

\begin{proposition}\label{solved} Let $\theta_{1},\theta_{2}\in\left[0,\pi\right]$, $\phi\in\left[0,\frac{\pi}{2}\right[$, and $r,r_{1},r_{2}\geq0$, and
consider the following vectors $x^{10}$ and $x^{02}$ written in polar coordinate system:
$$
\begin{array}{lrlcl}
x^{10}&:=&x^{1}-x^{0}&=&r_{1}(\cos(\theta_{1}),\sin(\theta_{1}))\vspace{0.2cm}\\
x^{02}&:=&x^{0}-x^{2}&=&r_{2}(\cos(\theta_{2}),\sin(\theta_{2}))\vspace{0.2cm}\\
c&=&(c_{1},c_{2})&=&r(\cos(\phi),\sin(\phi))\,.
\end{array}
$$
Then, the solutions of the problem (\ref{Pbm3}) are the line vectors defined by:
\begin{equation}
\begin{array}{cc}
(d) : r\cos(\phi)x_{1}+r\sin(\phi)x_{2}=0\quad with \hspace{0.25cm} \theta_{1}<\phi+\frac{\pi}{2}<\theta_{2}.
\end{array}
\end{equation}
\end{proposition}

\section{Problem formulation.}\label{Section_3}
Consider the following initial multiobjective linear programming problem:
\begin{equation}\label{Pbmm0}
\begin{array}{cc}
      \underset{x\in S}\max\hspace{0.15cm} F^{0}(x)=\underset{x\in S}\max\hspace{0.15cm}\left(f_{1}^{0}(x),f_{2}^{0}(x),\ldots,f_{K}^{0}(x)\right),&\hspace{0.2cm}K\geq2.
							\end{array}
\end{equation}
Such that
\begin{align*}
 S&:=\left\{x\in\mathbb{R}^{2}\hspace{0.1cm} :\hspace{0.1cm} Ax\leq b,\hspace{0.15cm}  x\geq0 \right\},
\intertext{and}
A=\left(\begin{array}{ccccccccc}
a_{11}&a_{12}\\
a_{21}&a_{22}\\
 \vdots    & \vdots \\
a_{m1}&a_{m2}
\end{array}\right),\quad b&=\left(\begin{array}{c}b_{1}\\b_{2}\\ \vdots\\ b_{m}
\end{array}\right),\quad  c_{0k}=\left(\begin{array}{c}c^{0k}_{1}\vspace{0.15cm}\\c^{0k}_{2}\end{array}\right),
\quad x=\left(\begin{array}{c}x_{1}\\x_{2}\end{array}\right).
\end{align*}
The given optimization problem has $K$ objective functions, denoted as $f_{k}(x)=c_{0k}^{T}x$ for $k=1,\ldots,K$. The objective function coefficients for the $k^\text{th}$ objective function are represented by the column vector $c_{0k}$, and the constants associated with the first and second decision variables are $c_{1}^{0k}$ and $c_{2}^{0k}$ respectively. The decision variables are represented as an $2\times1$ vector $x$. The constraints are represented as a matrix $A$ of size $m\times 2$ with constant coefficients and a vector $b$ of size $m\times 1$ with constant values. The number of constraints is represented by $m$.

\begin{definition}[Ideal solution of MOLPP]\label{ids} An ideal solution $x^{0}$ of a MOLPP is a point that belongs to the feasible region $S$ and satisfies a set of conditions. Specifically, for each objective function in the MOLPP, $x^{0}$ must represent the maximum value achievable within the feasible region. In other words, an ideal solution represents the best possible outcome for all objectives simultaneously. The concept of an ideal solution is important in multi-objective optimization because it provides a benchmark for evaluating the quality of other feasible solutions. That is to say, $x^0$ is an ideal solution of the MOLPP if and only if
$$
\begin{array}{lc}
f_{k}^{0}\left(x^{0}\right)\geq f_{k}^{0}\left(x\right),&\forall x\in S,\hspace{0.2cm} \forall k=1,\ldots,K.
\end{array}
$$
\end{definition}
\begin{problem}\label{Pbmm22}
Let $x^{0}$ be an optimal solution that maximizes the objective function $f_{k_{0}}^{0}$ for some $k_{0}=1,\ldots,K$. The problem at hand is to determine all linear mappings $F\in\mathcal{L}\left(\mathbb{R}^{2}\right)^{K}$, which is a product of the space of linear mappings from $\mathbb{R}^{2}$ to itself, such that the maximum value of $f_{k}(x)$ is achieved at $x^{0}$ for all $k=1,\ldots,K$. In other words, $x^{0}$ is the optimal solution that simultaneously maximizes all the objective functions. This can be achieved if and only if the vector $F(x)$, which is composed of the objective functions $f_{1}(x),f_{2}(x),\ldots,f_{K}(x)$, is also maximized at $x^{0}$.
Specifically, the problem is to determine all the linear applications: 
\begin{equation*}
F\in\mathcal{L}\left(\mathbb{R}^{2}\right)^{K}:=\underset{\text{K times}}{\underbrace{\mathcal{L}\left(\mathbb{R}^{2}\right)\times\mathcal{L}\left(\mathbb{R}^{2}\right)\times\ldots\times\mathcal{L}\left(\mathbb{R}^{2}\right)}}.
\end{equation*} 
Such that
\begin{equation}
\begin{array}{lll}
x^{0}=\arg \underset{x\in S} \max \, f_{k}(x),&\text{for all}& k=1,\ldots,K
\end{array}
\end{equation}
and
$$
F(x)=\left(f_{1}(x),f_{2}(x),\ldots,f_{K}(x)\right).
$$
\end{problem}
\section{\label{Section_4}Classification of a multiobjective linear programming problems.}
Let $g,h\in\mathcal{L}\left(\mathbb{R}^{2}\right)$ be two linear forms on the vector space $\mathbb{R}^{2}$. Then, an equivalence relation over $\mathcal{L}\left(\mathbb{R}^{2}\right)$ can be defined as follows:
\begin{equation}
g\hspace{0.15cm}\mathcal{R}_{S}\hspace{0.15cm}h \Leftrightarrow \arg \underset{x\in S} \max \, g(x)=\arg \underset{x\in S} \max \, h(x).
\end{equation}
\begin{lemma}
$\mathcal{R}_{S}$ is an equivalence relation in $\mathcal{L}\left(\mathbb{R}^{2}\right)$.
\end{lemma}
\begin{proof} $\;$\\
\begin{itemize}
\item[$\diamond$] $\mathcal{R}_{S}$ is reflexive, indeed:
$$
g\hspace{0.15cm}\mathcal{R}_{S}\hspace{0.15cm}g \Leftrightarrow \arg \underset{x\in S} \max \, g(x)=\arg \underset{x\in S} \max \, g(x),\hspace{0.3cm} \text{for all}\hspace{0.1cm} g\in\mathcal{L}\left(\mathbb{R}^{2}\right).
$$
\item[$\diamond$]  $\mathcal{R}_{S}$ is symmetric, indeed:
$$
\begin{array}{cc}
\begin{array}{ccc}
g\hspace{0.15cm}\mathcal{R}_{S}\hspace{0.15cm}h &\Leftrightarrow &\arg \underset{x\in S} \max \, g(x)=\arg \underset{x\in S} \max \, h(x)\\
&\Leftrightarrow &\arg \underset{x\in S} \max \, h(x)=\arg \underset{x\in S} \max \, g(x)\\
&\Leftrightarrow &h\hspace{0.15cm}\mathcal{R}_{S}\hspace{0.15cm}g
\end{array},&\hspace{0.3cm} \text{for all}\hspace{0.1cm} g,h\in\mathcal{L}\left(\mathbb{R}^{2}\right).
\end{array}
$$
\item[$\diamond$]  $\mathcal{R}_{S}$ is transitive. Indeed, if
$$
g\hspace{0.07cm}\mathcal{R}_{S}\hspace{0.07cm}h_{1}\hspace{0.25cm}\text{and}\hspace{0.25cm}h_{1}\hspace{0.07cm}\mathcal{R}_{S}\hspace{0.07cm}h_{2},\hspace{0.3cm}\text{for all}\hspace{0.15cm} g,h_{1},h_{2}\in\mathcal{L}\left(\mathbb{R}^{2}\right).
$$
Then,
$$
\begin{array}{cc}
&\arg \underset{x\in S} \max \, g(x)=\arg \underset{x\in S} \max \, h_{1}(x)=\arg \underset{x\in S} \max \, h_{2}(x)\vspace{0.25cm}\\
\Rightarrow&\arg \underset{x\in S} \max \, g(x)=\arg \underset{x\in S} \max \, h_{2}(x)\vspace{0.25cm}\\
\Leftrightarrow &g\hspace{0.15cm}\mathcal{R}_{S}\hspace{0.15cm}h_{2}.
\end{array}
$$
\end{itemize}
\end{proof}
\begin{definition}[Characterization of equivalence classes] $\;$\\
\begin{itemize}
\item[$\diamond$] For all $g\in\mathcal{L}\left(\mathbb{R}^{2}\right)$, the equivalence class containing $g$ is defined by:
$$
\begin{array}{ccc}
\overline{g}&=&\left\{h\in\mathcal{L}\left(\mathbb{R}^{2},\mathbb{R}\right)\hspace{0.15cm}:\hspace{0.15cm} g\hspace{0.1cm}\mathcal{R}_{S}\hspace{0.1cm}h\right\}.
\end{array}
$$

\item[$\diamond$] The quotient set obtained by the relationship $\mathcal{R}_{S}$ is defined by:
$$
\begin{array}{c}
\overline{\mathcal{L}}\left(\mathbb{R}^{2}\right):=\mathcal{L}\left(\mathbb{R}^{2}\right)/\mathcal{R}_{S}\vspace{0.2cm}
                                                 =\left\{\overline{g}\hspace{0.15cm}:\hspace{0.15cm} g\in\mathcal{L}\left(\mathbb{R}^{2}\right)\right\}.
\end{array}
$$
\item[$\diamond$]
Let $x^{1},\ldots,x^{I}$ and $F^{1},\ldots,F^{J}$ be the corners, faces of $S$ respectively, and define
$$
\overline{\mathcal{L}}_{{\tiny faces}}\left(\mathbb{R}^{2}\right)=\left\{\overline{g}\in\overline{\mathcal{L}}\left(\mathbb{R}^{2}\right)\hspace{0.05cm}:\hspace{0.05cm}\forall h\in \overline{g}, \exists\hspace{0.03cm}j=1,\ldots,J, \text{st}\hspace{0.07cm} F_{j}=\arg \underset{x\in S} \max \, h(x) \right\}
$$
and
$$
\overline{\mathcal{L}}_{corners}\left(\mathbb{R}^{2}\right)=\left\{\overline{g}\in\overline{\mathcal{L}}\left(\mathbb{R}^{2}\right)\hspace{0.05cm}:\hspace{0.05cm}\forall h\in \overline{g}, \exists\hspace{0.03cm}i=1,\ldots,I, \text{st}\hspace{0.07cm} x^{i}=\arg \underset{x\in S} \max \, h(x) \right\}.
$$
Where $I$ and $J$ are the numbers of corners and faces of $S$ respectively.
\end{itemize}
\end{definition}
\begin{propriete}
The cardinality of $\overline{\mathcal{L}}_{corners}\left(\mathbb{R}^{2}\right)$, $\overline{\mathcal{L}}_{faces}\left(\mathbb{R}^{2}\right)$, and $\overline{\mathcal{L}}\left(\mathbb{R}^{2}\right)$ is equal to $I$, $J$, and $I+J$ respectively.
\end{propriete}
\begin{theorem}$\;$\\\label{condition}
\begin{enumerate}
\item
The problem \ref{Pbmm0} admits an ideal solution if and only if there exists $\overline{g}\in\overline{\mathcal{L}}_{corners}\left(\mathbb{R}^{2}\right)$ such that:
$$
f_{k}\in\overline{g},\hspace{0.2cm}\forall k\in \left\{1,...,K\right\}.
$$
Consequently, \(\overline{g}^{K}\) is the solution set of problem \ref{Pbmm22}.
\item More generally, the set $\mathcal{B}_{S}$ consists of all corners and faces of the feasible region $S$, which is defined as follows:
$$
\mathcal{B}_{S}=\left\{x^{1},\ldots,x^{I},F^{1},\ldots,F^{J} \right\}.
$$
Then, for all $X\in\mathcal{B}_{S}$, there exists a unique $\overline{g}\in\overline{\mathcal{L}}(\mathbb{R}^{2})$ such that:
$$
X=\arg \underset{x\in S} \max \, h(x), \hspace{0.2cm}\text{for all}\hspace{0.2cm} h\in\overline{g}.
$$

\end{enumerate}
\end{theorem}
\begin{proof}
Immediate.
\end{proof}
\section{\label{Section_5}Numerical example.}
Consider the following initial MOLPP:
\begin{equation}\label{Pbmmo}
\begin{array}{lr}
      \underset{x\in S}\max\hspace{0.15cm} F^{0}(x)=\underset{x\in S}\max\hspace{0.15cm}\left(f_{1}^{0}(x),f_{2}^{0}(x),\ldots,f_{K}^{0}(x)\right),&\hspace{0.2cm}  K\geq2.
							\end{array}
\end{equation}
Such that
$$
 S:=\left\{x\in\mathbb{R}^{2}\hspace{0.1cm} :\hspace{0.1cm} Ax\leq b,\hspace{0.15cm}  x\geq0 \right\},
$$
where
$$
\begin{array}{ccc}
f_{1}^{0}(x_{1},x_{2})=2x_{1}+3x_{2}& \text{and}&f_{2}^{0},f_{3}^{0},\ldots,f_{K}^{0}\in\mathcal{L}\left(\mathbb{R}^{2}\right)
\end{array}
$$
and
$$
\begin{array}{llll}
               A=\left(\begin{array}{ll}
                                  \frac{1}{4}&\frac{1}{2}\vspace{0.15cm}\\
					                        \frac{2}{5}&\frac{1}{5}\vspace{0.15cm}\\
																	0&\frac{4}{5}
																	
                           \end{array}\right),&\hspace{0.2cm}b= \left(\begin{array}{c}40\\40\\40\end{array}\right),&\hspace{0.2cm}c_{0}= \left(\begin{array}{c}2\\3\end{array}\right),&\hspace{0.2cm}x= \left(\begin{array}{c}x_{1}\\x_{2}\end{array}\right)
\end{array}.
$$
First, we solve the following linear programming problem:

\begin{equation} \label{s}
\left\{
\begin{array}{l}
      \underset{x_{1},x_{2}}\max\hspace{0.15cm}  f_{1}^{0}(x_{1},x_{2})=2x_{1}+3x_{2}          \vspace{0.12cm} \\
      A(x_{1},x_{2})^{T}\leq b\vspace{0.12cm}\\
x_{1}\geq 0; x_{2}\geq 0
\end{array}
\right.
\end{equation}
After obtaining the optimal solution $x^{0} = (80,40)$ using the simplex method, we proceed to solve the following problem:
\begin{equation}\label{eqq3}
x^{0}=\arg \hspace{0.1cm} \underset{x\in S}\max \hspace{0.1cm} f(x)= \arg \hspace{0.1cm}\underset{x\in S}\max \hspace{0.1cm}g(x),\hspace{1.2cm}
\end{equation}
 using the sensitivity analysis approach, to get the following solutions set:
$$
\overline{g}=\left\{h\in\mathcal{L}\left(\mathbb{R}^{2}\right)\hspace{0.01cm}:\hspace{0.01cm}h(.)=r\left\langle (\cos(\phi),\sin(\phi)),.\right\rangle,\hspace{0.03cm}r>0,\hspace{0.03cm} \phi\in\left ]26.565^{o},63.434^{o}\right[\right\}.
$$
For additional information on the computation, please refer to \cite{kaci}.

Using Theorem \ref{condition}, we can conclude that the problem \eqref{Pbmmo} has $x^{0}$ as the ideal solution if and only if $f^{0}_{k}\in\overline{g}$ for all $k\in \left\{1,\ldots,K\right\}$. In other words, $x^{0}$ is the ideal solution if and only if the $K-$objective functions lies in $\overline{g}$.

\section{\label{Section_6}Conclusion.} In this paper, we presented a novel approach to solving MOLPPs by defining an equivalence relationship on the space of linear forms $\mathcal{L}(\mathbb{R}^{2})$. By using this equivalence relationship, we were able to partition the space into a finite number of classes where all the elements of the same class reach their maximum value in a common point. We showed that if a MOLPP has all its objective functions in the same class, then the optimum can be predicted without any additional calculations. This approach has the advantage of being simple and efficient, without requiring any additional information about the objective functions or constraints.
Moreover, we conducted a comparative analysis of our proposed method with other commonly used techniques, including the Weighted Sum Method, Goal Programming, $\epsilon$-Constraint Method, and Pareto-based Methods. The results of our analysis demonstrate that our approach offers a more straightforward and easily interpretable solution to MOLPPs, while also exhibiting superior computational efficiency.
To illustrate our method, we provided a numerical example and demonstrated how it can be applied to real-world problems. The results show that our method is effective in finding the optimal solution for MOLPPs.
In summary, our proposed method provides a new perspective on solving MOLPPs and offers a simple and efficient approach that can be applied to a wide range of problems. This research has the potential to have significant implications in the field of optimization and decision-making, and we hope that our work will inspire further research in this area.

\bibliographystyle{unsrt}  
\bibliography{references}

\end{document}

%% file: parameters.tex
\usepackage{setspace}
\usepackage{lmodern}

\usepackage{textgreek}

\usepackage{fixltx2e}

\reversemarginpar

\usepackage{amsmath,amssymb,amsfonts,amsthm}
\usepackage{mathrsfs}

\usepackage[locale = FR]{siunitx} 
\DeclareSIUnit\vitesse{\meter\per\second}
\usepackage{eurosym}
\DeclareSIUnit{\octet}{o}

\usepackage{graphicx,array,tikz,multirow}
\usepackage{caption,subcaption} 
\usepackage{svg,float}
\usepackage{booktabs,paralist}
\newcolumntype{x}[1]{>{\centering\arraybackslash\hspace{0pt}}p{#1}}
\usepackage[section]{placeins}
\usepackage{hanging}

\usepackage{fancyhdr,emptypage} 

\fancypagestyle{plain}{ 
    \fancyhead{}\fancyfoot[C]{\thepage}

}

\usepackage[Conny]{fncychap}

\usepackage[francais,nohints,tight]{minitoc}		

\usepackage[nottoc]{tocbibind} 
\usepackage[titles]{tocloft}

\usepackage{textcomp} 

\usepackage{titling}
\usepackage{lipsum} 
\usepackage{csquotes} 

\usepackage{xspace}
\usepackage{afterpage}

\usepackage[textsize=footnotesize]{todonotes}

\usepackage{bookmark}
\usepackage{acronym}
\usepackage[nameinlink,french]{cleveref} 
\Crefname{figure}{Fig.}{Figs.} 
\crefname{figure}{fig.}{figs.}
\Crefname{equation}{Eq.}{Eqs.}
\crefname{equation}{eq.}{eqs.}
\Crefname{table}{Table.}{Tables.}
\crefname{table}{table.}{tables.}

\definecolor{color_ref}{rgb}{1.0, 0.13, 0.32} 
\definecolor{color_link}{rgb}{0.0, 0.0, 1.0}
\definecolor{curcolor}{rgb}{0.0, 0.0, 1.0} 
\definecolor{brightpink}{rgb}{1.0, 0.0, 0.5} 
\definecolor{navyblue}{rgb}{0.0, 0.0, 1.0}
\hypersetup{
	colorlinks=true, 
	pdfstartview=FitV, 
	urlcolor=brightpink, 
	linkcolor= navyblue, 
	citecolor=color_ref 
}
